\documentclass{amsart}
\usepackage{enumerate}
\usepackage[latin1]{inputenc}

\newtheorem{e-proposition}[theorem]{Proposition}

\newtheorem{e-definition}[theorem]{Definition\rm}
\newtheorem{remark}{\it Remark\/}

\def\og{\leavevmode\raise.3ex\hbox{$\scriptscriptstyle\langle\!\langle$~}}
\def\fg{\leavevmode\raise.3ex\hbox{~$\!\scriptscriptstyle\,\rangle\!\rangle$}}

\def\ds{\displaystyle}
\def\R{\mathbb{R}}
\def\be{\begin{equation}}
\def\ee{\end{equation}}

\def\n{{\boldsymbol n}}
\def\p{\partial}
\def\grad{\boldsymbol{\nabla}}
\def\div{\grad\cdot}

\def\Om{\Omega}

\def\a{\alpha}

\def\x{{\boldsymbol x}}
\def\d{{\rm d}}
\def\ov#1{\overline{#1}}
\def\un#1{\underline{#1}}
\def\wh#1{\widehat{#1}}

\def\Ee{\mathcal{E}}
\def\Zz{\mathcal{Z}}
\def\Pp{\mathcal{P}}
\def\Dd{\mathcal{D}}
\def\Mm{\mathcal{M}}

\def\s{{\boldsymbol{s}}}
\def\e{{\boldsymbol{e}}}
\def\h{{\boldsymbol{h}}}

\def\V{{\bf V}}
\def\L{{\boldsymbol \Lambda}}
\def\bv{{\boldsymbol v}}
\def\bu{{\boldsymbol u}}
\def\bp{{\boldsymbol p}}
\def\0{{\bf 0}}

\begin{document}
\centerline{}


\title[Gradient structure for two-phase porous media flows]{The gradient flow structure for incompressible immiscible two-phase flows in porous media}


\author[C. Canc\`es]{Cl\'ement Canc\`es}
\address{
Cl\'ement Canc\`es ({\tt cances@ljll.math.upmc.fr})
\begin{enumerate}
\item Sorbonne Universit\'es, UPMC Univ Paris 06, UMR 7598, Laboratoire Jacques-Louis Lions, F-75005, Paris, France 
\item CNRS, UMR 7598, Laboratoire Jacques-Louis Lions, F-75005, Paris, France
\end{enumerate}
}

\author[T. O. Gallou\"et]{Thomas O. Gallou\"et}
\address{
Thomas O. Gallou\"et ({\tt thomas.gallouet@inria.fr})
\begin{enumerate}
\item Universit\'e Libre de Bruxelles (ULB), Brussels, Belgium
\item Project-Team MEPHYSTO, Inria Lille - Nord Europe, Villeneuve d'Ascq, France
\end{enumerate}
}

\author[L. Monsaingeon]{L\'eonard Monsaingeon}
\address{
L\'eonard Monsaingeon ({\tt leonard.monsaingeon@ist.utl.pt})
\begin{enumerate}
\item 
CAMGSD, Instituto Superior T\'ecnico, Universidade de Lisboa, Av. Rovisco Pais, 1049-001 Lisboa, Portugal
\end{enumerate}
}


\begin{abstract}

We show that the widely used model governing the motion of two 
incompressible immiscible fluids in a possibly heterogeneous porous 
medium has a formal gradient flow structure. More precisely, the 
fluid composition is governed by the gradient flow of some non-smooth 
energy. Starting from this energy together with a dissipation potential, we recover the 
celebrated Darcy-Muskat law and the capillary pressure law governing 
the flow thanks to the principle of least action. Our interpretation does not 
require the introduction of any algebraic transformation like, e.g., the global 
pressure or the Kirchhoff transform, and can be transposed to the case of more phases. 

\end{abstract}

\maketitle


\section{Introduction}
\label{}

\subsection{General motivations}
The models for multiphase porous media flows have been widely studied in 
the last decades since they are of great interest in several fields of applications, 
like e.g. oil-engineering, carbon dioxide sequestration, or nuclear waste repository 
management. We refer to the monographs~\cite{Bear72,BB90} for an extensive 
discussion on the derivation of models for porous media flows, and to \cite{AS79,CJ86,AKM90,CHM06} 
for numerical and mathematical studies.
\medskip

More recently, F. Otto showed in his seminal work~\cite{Otto01} that the so-called \emph{porous 
medium equation}: 
$$
\p_t \rho - \Delta \rho^m = 0 \quad \text{ for } (\x,t) \in \R^N \times \R_+ \text{ and } m >1, 
$$ 
which is a very simplified model  corresponding to the case of an isentropic gas flowing  within a porous 
medium, can be reinterpreted in a physically relevant way as the gradient flow of 
the free energy with respect to some Wasserstein metric in the space of Borel probability measures. Extensions to more general degenerate 
parabolic equations were then proposed for example in~\cite{Agueh05,Lisini09}. 
See also for instance ~\cite{Bla14_KS} or~\cite{KMX15} for the interpretation of 
some dissipative systems as gradient flows in Wasserstein metrics. 
\medskip

In this note, we will focus on the model governing the motion of 
an incompressible immiscible two-phase flow in a possibly heterogeneous porous medium, 
that will appear in the sequel as~\eqref{eq:cons} and \eqref{eq:Darcy-Muskat}--\eqref{eq:pore-sat}.
This model is relevant for instance for describing the flow of oil and water, whence the 
subscripts $o$ and $w$ appearing in the sequel of this note, within a rock that is possibly made of 
several rock-types. Our goal is to show that, at least formally, 
this model can be reinterpreted as the gradient flow of some singular energy. This will motivate the use of structure-preserving numerical 
methods inspired from~\cite{CG_VAGNL} to this model in the future.

\medskip

Our approach is inspired from the one of A. Mielke~\cite{Mie11} and, more closely, 
to the one of M.~A. Peletier~\cite{Pel-lecture}. The basic recipe for variational modeling 
is recalled in~\S\ref{ssec:recipe}, then its ingredients are identified in \S\ref{sec:2}. 
This approach is purely formal, but it can be made 
rigorous under some unphysical strict positivity assumption on the phase 
mobilities $\eta_o, \eta_w$ defined below. We will remain sloppy about regularity issues all along this 
note. 

\subsection{The recipe of getting formal variational models }\label{ssec:recipe}

Here we recall very briefly the main ingredients needed for defining a 
formal gradient flow. 
\begin{enumerate}[{\bf i.}]
\item The {\em state space} $\Mm$ is the set where the solution of the gradient flow can evolve.
\item At a point $\s \in \Mm$, the tangent space $T_\s\Mm$, to whom would belong $\p_t \s$, 
is identified in a non-unique way with a so-called \emph{process space} $\Zz_\s$ (that might depend on $\s$).
More precisely, we assume that for each $\s\in\Mm$ there exists an onto linear application $\Pp(\s):\Zz_\s \to T_\s\Mm$.
\item The {\em energy functional} $\Ee: \Mm \to \R\cup\{+\infty\}$ admits a (local) sub-differential $\p_\s \Ee(\s) \subset \left(T_\s\Mm\right)^\ast$ at 
 $\s \in \Mm$.
\item The {\em dissipation potential} $\Dd$ is such that, for all $\s \in \Mm$ and all $\V \in \Zz_\s$, one has 
$\Dd(\s;\V) \ge 0$. It is supposed to be convex and coercive w.r.t. to its second variable.
\item The initial data $\s^0$ belongs to $\Mm$.
\end{enumerate}
All these ingredient being defined, we obtain from the {\em principle of least action} that 
$\s : \R_+ \to \Mm$ is the gradient flow of the energy $\Ee$ for the dissipation $\Dd$ if 
\begin{subequations}
\be
\p_t \s = \Pp(\s) \V
\ee 
where  
\be\label{eq:leastaction}
\V\in \underset{\wh \V \in \Zz_{\s} }{\rm argmin}
\left( \max_{\h \in \p_\s \Ee(\s)} \left(
\Dd\big(\s(t);\wh \V(t)\big)  + \Big\langle \h \, ,  \, \Pp(\s) \wh \V \Big\rangle_{(T_\s\Mm)^\ast, T_\s\Mm} \right)\right).
\ee
\end{subequations}
The formula~\eqref{eq:leastaction} means that a gradient flow is lazy and smart: 
the motion aims to minimize the dissipation while maximizing the decay of the energy.
We refer to~\cite{Mie11,Pel-lecture} for additional material on such a formal modeling and to~\cite{AGS08} 
for an extensive (and rigorous) discussion on gradient flows in metric spaces.

\section{Variational modeling for two-phase flows in porous media}\label{sec:2}

\subsection{State space and process space}\label{ssec:state-process}

Let $\Om$ be an open subset of $\R^N$ representing a (possibly heterogeneous) \emph{porous medium}, 
let $\phi : \Om \to (0,1)$ be a measurable function (called \emph{porosity}) 
such that $\un \phi \le \phi(\x) \le \ov \phi$ for a.e. $\x \in \Om$
for some constants $\un \phi, \ov \phi \in (0,1)$, 
and let $\un s_{o}, \un s_{w}:\Om \to [0,1)$ be two measurable functions 
(so-called \emph{residual saturations}) such that 
$\un s_{o}(\x) + \un s_{w}(\x) < 1$ for a.e. $\x \in \Om.$
In what follows, we denote by 
$$
\ov s_o(\x) = 1-\un s_w(\x), \qquad \ov s_w(\x) = 1 - \un s_o(\x), \qquad \text{ for a.e. } \x \in \Om.
$$
For almost all $\x \in \Om$, we denote by 
$$
\Delta_\x = \Big\{ \s=(s_o,s_w) \in \R^2 \left| \; s_o + s_w = 1 \text{ with }
\un s_{\a}(\x) \le  s_\a \le \ov s_\a(\x) \text{ for } \a \in \{o,w\}\Big\}.\right. 
$$

Let $\s^0 = (s_o^0, s_w^0)$ be a given initial saturation profile, we denote by $m_\a$ ($\a \in \{o,w\}$) 
the volume occupied by the phase $\a$ in the porous medium, i.e., 
$$
m_o = \int_\Om \phi(\x) s_o^0(\x) \d\x, \quad \text{ and }\quad  m_w = \int_\Om \phi(\x) s_w^0(\x) \d\x.
$$ 
For simplicity, we restrict our attention to the case where the volume of each phase is preserved: 
no source term and no-flux boundary conditions (otherwise, non-autonomous gradient flows should be considered). 
Hence the saturation profile lies at each time in the so-called state space $\Mm$, defined by 
$$
\Mm = \left\{ \s = (s_o,s_w)  \; \left| \; s_\a : \Om \to \R_+ \; \text{ with } \; 
\int_\Om \phi(\x) s_\a(\x) \d\x = m_\a\;  \text{ for } \a \in \{o,w\}  \right\}.\right.
$$
Let us now describe the processes that allow to transform the saturation profile. 
We denote by 
$$
\Zz_\s = \Big\{ \V = (\bv_o, \bv_w) \; \left| \; \bv_\a: \Om \to \R^N\; \text{ with } \bv_\a \cdot \n =0 \text{ on } \p\Om \Big\}\right.
$$
the \emph{process space}  of the admissible 
processes for modifying a saturation profile $\s\in\Mm$. The identification between $\V =(\bv_o, \bv_m)  \in \Zz_\s$ and 
$ \boldsymbol{\dot\s} = (\dot s_o,\dot s_w) \in T_\s\Mm$ is made through the onto operator $\Pp(\s):\Zz_\s \to T_\s\Mm$ defined 
by 
\be\label{eq:P(s)}
\Pp(\s)\V = \left(-\frac1{\phi}\div \bv_o \, ; \, -\frac1{\phi}\div \bv_w \right), \qquad \forall \V \in \Zz_\s.
\ee

Since $\p_t\s \in T_\s\Mm$, the relation~\eqref{eq:P(s)} yields the existence of some \emph{phase filtration speeds} $(\bv_o, \bv_w) \in \Zz_\s$ 
such that the following {\em continuity equations}  hold: 
\be\label{eq:cons}
\phi \p_t s_\a + \div\bv_\a = 0, \qquad \a \in \{o,w\}.
\ee
The relation~\eqref{eq:cons} must be understood as the local volume conservation of each phase $\a \in \{o,w\}$.
Finally, the duality bracket $\langle \cdot \, ,  \, \cdot \rangle_{(T_\s\Mm)^\ast, T_\s\Mm}$ is 
given by
\begin{align*}
\langle \h ,  \dot\s \rangle_{(T_\s\Mm)^\ast, T_\s\Mm} =& \sum_{\a \in \{o,w\}} \int_\Om \phi h_\a \dot s_\a \\
=& - \sum_{\a \in \{o,w\}} \int_\Om h_\a \div \bv_\a =  \sum_{\a \in \{o,w\}} \int_\Om \grad h_\a \cdot \bv_\a.
\end{align*}
\subsection{About the energy}\label{ssec:nrj}

For a.e. $\x \in \Om$, we assume  $\pi(\cdot,\x): [\un s_{o}(\x),\ov s_o(\x)] \to \R$ 
to be a maximal monotone graph whose restriction 
$\pi_{|_{(\un s_o,\ov s_o)}}(\cdot, \x)$ to the open interval $(\un s_o(\x),\ov s_o(\x))$ 
is an increasing (single-valued) function belonging to $L^1(\un s_o(\x),\ov s_o(\x))$.  
In particular, $\pi^{-1}(\cdot,\x) : \R \to [\un s_{o}(\x),\ov s_o(\x)]$ is a single valued function.

We denote by $\Pi:\R  \times \Om \to \R\cup\{+\infty\}$ the (strictly convex w.r.t. its first variable) function defined by 
$$
\Pi(s_o,\x) = \begin{cases}
\ds \int^{s_o}_{\sigma(\x)} \pi(a,\x)\d a  - (\rho_o - \rho_w) s g z  & \text{ if } s_o \in [\un s_{o}(\x),\ov s_{o}(\x)], \\
+ \infty & \text{ otherwise}, 
\end{cases}
$$
where, denoting by $\e_z$ the downward unit normal vector of $\R^N$, we have set $z = \x \cdot \e_z$, 
and where $g$ and $\rho_\a$ denote the gravity constant and the density of the phase $\a$ respectively, 
and where $\sigma$ is such that $\x \mapsto \pi(\sigma(\x),\x) - (\rho_o - \rho_w) g z$ is constant.
Since $\pi_{|_{(\un s_o,\ov s_o)}}(\cdot, \x) \in L^1(\un s_o(\x),\ov s_o(\x))$, we get that  $\Pi(\un s_{o}(\x),\x)$ and $\Pi(\ov s_{o}(\x), \x)$ are finite for a.e. $\x \in \Om$.

The \emph{volume energy} function $E: \R^2 \times \Om \to \R \cup\{+\infty\}$ 
is defined by
\be\label{eq:E}
E(\s,\x) = \begin{cases}
\Pi(s_o,\x) & \text{ if } \s = (s_o,s_w)\in \Delta_\x, \\
+\infty & \text{ otherwise}.
\end{cases}
\ee
The function $E(\cdot, \x)$ is convex and finite on $\Delta_\x$ for a.e. $\x \in \Om$. Its sub-differential is given by 
$$
\p_\s 
E(\s,\x) = 
\begin{cases}
\Big\{ (h_o ,h_w) \in \R^2\; \left| \; h_o - h_w + (\rho_o - \rho_w )g z \in \pi(s_o,\x) \Big\}\right. & \text{ if } \s  \in \Delta_\x, \\
\emptyset & \text{ otherwise}.
\end{cases}
$$

Finally, we can define the so-called \emph{global energy} $\Ee : \Mm \to \R \cup  \{+\infty\}$ by 
\be\label{eq:Ee}
\Ee(\s) = \int_\Om \phi(\x) E(\s(\x), \x) \d \x, \qquad \forall \s=(s_o,s_w) \in \Mm. 
\ee
The saturation profile $\s \in \Mm$ is of finite energy $\Ee(\s) < \infty$ if and only if $\s(\x) \in \Delta_\x$ for a.e. $\x \in \Om$. For $\s \in \Mm$ with finite energy one can check that the local sub-differential $\p_\s \Ee(\s)$ of $\Ee$ at $\s$ is given by 
\begin{multline}
\label{eq:pEe}
\p_\s \Ee(\s) = \Big \{\h = (h_o,h_w) : \Om \to \R^2  \text{ s.t. }  \\
 h_o - h_w + (\rho_o - \rho_w )g z \in \pi(s_o,\x)\; \text{ for a.e. } \x \in \Om \Big\}.
\end{multline}

\subsection{About the dissipation}\label{ssec:dissip}

The {\em permeability tensor} field $\L\in L^\infty(\Om;\R^{N\times N})$ is assumed to be such that 
$\L(\x)$ is a symmetric and positive matrix for a.e. $\x \in \Om$.
Moreover, we assume that there exist $\lambda_\star, \lambda^\star \in \R_+^\ast$  such that 
$$
\lambda_\star | \bu |^2 \le \L(\x) \bu \cdot \bu \le \lambda^\star | \bu |^2, \qquad \text{for all } 
\bu \in \R^N \text{ and a.e. } \x \in \Om. 
$$
This ensures that $\L(\x)$ is invertible for a.e. $\x \in \Om$. Its inverse is denoted by $\L^{-1}(\x)$.

We also need the two Carath\'eodory functions $\eta_o, \eta_w : \R \times \Om \to \R_+$ 
--- the so-called {\em phase mobilities} ---  such that  $\eta_\a(\cdot, \x)$ are
Lipschitz continuous and nondecreasing on $\R_+$ for a.e. $\x \in \Om$ and $\a \in \{o,w\}$.
Moreover, we assume that 
$\eta_\a(s,\x) = 0$ if $s \le \un s_\a(\x)$ and that  $\eta_\a(s,\x) >0$  if  $s > \un s_\a(\x)$.

Given $\s = (s_o,s_w) \in \Mm$ and $\V = (\bv_o,\bv_w)\in \Zz_\s$, 
we define the \emph{dissipation potential } $\Dd$ by 
$$
\Dd(\s,\V) = \frac12 \sum_{\a \in \{o,w\}} {\int_\Om} \frac{\L^{-1} \bv_\a \cdot \bv_\a}{\eta_\a(s_\a)} \d\x, \qquad 
\forall \s \in \Mm, \; \forall \V \in \Zz_\s.
$$
It is easy to check that dissipation is finite, i.e., $\Dd(\s,\V) < \infty$, iff 
$\bv_\a = \0  \text{ a.e. on } \{ \x \in \Om \; | \; s_\a(\x) \le \un s_\a(\x) \}.$

\subsection{Principle of least action and resulting equations}
Let us consider the gradient flow governed by the energy $\Ee$, 
the continuity equation~\eqref{eq:cons}, and the dissipation $\Dd$. 
Let $\s \in \Mm$ be a finite energy saturation profile, then 
because of the {\em principle of least action}~\eqref{eq:leastaction} and of the definition~\eqref{eq:P(s)} 
of the operator $\Pp(\s): \Zz_\s \to T_\s \Mm$, the process $\V=(\bv_o,\bv_w) \in \Zz_\s$ and the 
\emph{hydrostatic phase pressures} $\h = (h_o, h_w)$ must be chosen so that $(\V,\h)$ is the $\min-\max$ 
saddle-point of the functional 
\be\label{eq:legendre}
(\wh \V, \wh \h) \mapsto  \Dd(\s, \wh\V) - \sum_{\a \in \{o,w\}}\int_\Om \wh h_\a \div \wh \bv_\a \d\x.
\ee

One can first fix $\wh \h\in \p_\s \Ee(\s)$ and minimize w.r.t. $\V$. 
This provides 
\be
\label{eq:min_V_wrt_h}
 \underset{\substack{\wh\V \in \Zz}}{\rm argmin} 
\left( \Dd(\s, \wh\V) -\!\!\! \sum_{\a \in \{o,w\}}\int_\Om \wh h_\a \div \wh \bv_\a \d\x \right) 
= \left(- \eta_o (s_o) \L \grad \wh h_o, - \eta_w (s_w) \L\grad \wh h_w \right).
\ee
Injecting this expression in~\eqref{eq:legendre} and  maximizing w.r.t. $\wh \h\in \p_\s \Ee(\s)$, that is
minimizing
\be
\label{eq:min_subgradient_Wasserstein}
\h=\underset{\substack{\wh\h \in \p_\s\Ee(\s)}}{\rm argmin} \left(\frac{1}{2}\int_\Om\eta_\a(s_\a)\L \nabla \hat h_\a \cdot \nabla  \hat h_\a\right)
\ee
among all elements $\wh\h$ in the subdifferential $\p_\s\Ee(\s)$, yields 
\be\label{eq:total-flow}
- \div \Big( \bv_o + \bv_w \Big) = 0, \qquad  \bv_\a = - \eta_\a(s_\a) \L \grad h_\a.
\ee
In \eqref{eq:total-flow} the first condition follows from the constraint $\hat h \in\p_\s\Ee(\s)$ in \eqref{eq:min_subgradient_Wasserstein}, and the second one from \eqref{eq:min_V_wrt_h}.

Define the {\em phase pressures} $\bp = (p_o, p_w)$ by 
$
p_\a(\x) = h_\a(\x) + \rho_\a g z,
$
for a.e. $\x \in \Om$ and $\a \in \{o,w\}$, then
we recover the classical \emph{Darcy-Muskat law}: 
\be\label{eq:Darcy-Muskat} 
\bv_\a = - \eta_\a(s_\a) \L \grad\left( p_\a - \rho_\a g z \right), \quad \a \in \{o,w\}.
\ee
Moreover, it follows from~\eqref{eq:pEe} that the following \emph{capillary pressure relation} holds:
\be\label{eq:pc}
p_o(\x) - p_w(\x) \in \pi(s_o(\x),\x) \quad \text{ a.e. in } \Om.
\ee
We recover here the multivalued capillary pressure relation proposed in~\cite{CGP09,CP12}. 

Combining~\eqref{eq:cons} and~\eqref{eq:total-flow} easily gives $\p_t(s_o+s_w)=0$, so that the condition 
\be\label{eq:pore-sat}
s_o + s_w = 1 \quad \text{a.e. in }\Om, 
\ee
is preserved along time and the whole pore volume remains saturated by the two fluids.

Gathering~\eqref{eq:cons}, \eqref{eq:Darcy-Muskat}, \eqref{eq:pc} and~\eqref{eq:pore-sat} 
gives the usual system of equations governing immiscible incompressible two-phase flows in 
porous media~\cite{Bear72,CJ86,AKM90,Chen01,CP12}. \\

\begin{remark}
By similarity with the classical Wasserstein distance used in optimal mass transport \cite{Otto01} one could here 
endow the tangent space $T_\s\Mm$ at $\s \in \Mm$ with a weighted $\dot H^{-1}$-scalar product 
$$
\big(\dot \s_1, \dot \s_2 \big)_{T_\s\Mm} = \sum_{\a \in \{o,w\}}\int_\Om
 \eta_\a(s_\a) \L \grad h_{1,\a} \cdot\grad h_{2,\a}
  \d\x, \qquad 
$$
where, for $i \in \{1,2\}$ and $\a \in \{o,w\}$, we have set $\dot \s_i = (\dot s_{i,o}, \dot s_{i,w})$ and where $h_{i,\a}$ solves 
$$
- \div \left( \eta_\a(s_\a) \L \grad h_{i,\a} \right) = \dot s_{i,\a} \text{ in }\Om, \qquad 
\eta_\a(s_\a) \L \grad h_{i,\a}\cdot \n = 0 \text{ on }\p\Om. 
$$
Under some conditions on the functions $\eta_\a$ (see~\cite{DNS09}), this should allow us to 
consider $\Mm$ as a metric space endowed with the corresponding distance, but $\Ee$ is not locally $\lambda$-convex for this Riemannian structure. The minimization \eqref{eq:min_subgradient_Wasserstein} then consists in the selection of the subgradient with minimal norm.
\end{remark}

\section*{Acknowledgements}
This work was supported by the French National Research Agency ANR through grant ANR-13-JS01-0007-01
(Geopor project). TG acknowledges financial support from the European Research Council under the European
Community's Seventh Framework Programme (FP7/2014-2019 Grant Agreement QUANTHOM
335410). 
LM was supported by the Portuguese Science Fundation through FCT fellowship SFRH/BPD/88207/2012.

\end{document}